\documentclass[12pt]{article}
\usepackage{amsmath, amsthm, amssymb}

\title{Generic families of matrix pencils and
their bifurcation diagrams\footnotetext{This is the authors' version of a work that was published in Linear Algebra Appl. 332--334 (2001) 165--179.}}

\author{M$^{\mbox{\b{a}}}$
  Isabel Garc\'\i a-Planas\\
Dept. de Matem\`atica Aplicada I\\ Universitat
Polit\`ecnica de Catalunya\\ Miner{\'\i}a 1-1-3,
08038 Barcelona, Spain\\ igarcia@ma1.upc.es\\
\\ Vladimir V. Sergeichuk \\
        Institute of Mathematics,
        Tereshchenkivska St. 3 \\
        Kiev, Ukraine \\
        sergeich@ukrpack.net}
\date{}
\begin{document}
\maketitle

\begin{abstract}
V. I. Arnold (``On matrices depending on
parameters", {\it Russian Math. Surveys} 26, no.
2, 1971, 29--43) constructed smooth generic
families of matrices with respect to similarity
transformations depending smoothly on the entries
of matrices and got bifurcation diagrams of such
families with a small number of parameters. We
extend these results to pencils of matrices.
\end{abstract}

\def\newpic#1{%
   \def\emline##1##2##3##4##5##6{%
      \put(##1,##2){\special{em:point #1##3}}%
      \put(##4,##5){\special{em:point #1##6}}%
      \special{em:line #1##3,#1##6}}}
\newpic{}

\newcommand{\matr}[4]%
{\left[\genfrac{}{}{0pt}{}{#1}{#3}\,
\genfrac{}{}{0pt}{}{#2}{#4}\right]}

\renewcommand{\le}{\leqslant}
\renewcommand{\ge}{\geqslant}

\newcommand{\gl}{\mathop{\rm Gl}\nolimits}
\newcommand{\diag}{\mathop{\rm diag}\nolimits}

\newtheorem{theorem}{Theorem}[section]
\newtheorem{lemma}{Lemma}[section]
\newtheorem{corollary}{Corollary}[section]

\theoremstyle{remark}
\newtheorem{remark}{Remark}[section]
\newtheorem{example}{Example}[section]

\section{Introduction}
V. I. Arnold \cite{arn1} (see also \cite{arn2})
obtained a miniversal deformation of a Jordan
matrix; that is, a simplest possible canonical
form, to which not only a given square matrix
$A$, but also an arbitrary family of matrices
close to $A$, can be reduced by means of a
similarity transformation that depends smoothly
on the entries of $A$. Using this miniversal
deformation, Arnold (\cite{arn1}, \cite{arn2})
constructed bifurcation diagrams for generic
smooth one-, two- and three-parameter families of
matrices $A(\alpha_1,\dots,\alpha_n),\ n\le 3$;
that is, he described all possible types of
Jordan forms of $A(\alpha_1,\dots,\alpha_n)$ in a
neighborhood of $\vec 0$. The results are
important for applications in which one has
matrices that arise from physical measurements,
which means that their entries are known only
approximately.

Miniversal deformations of matrix pencils were
obtained in \cite{kag} and \cite{gar_ser}. In
this article we construct bifurcation diagrams
for generic smooth zero-, one- and two-parameter
families of pencils.

The case of zero parameter families is trivial.
(Of course, ``a zero parameter family of pencils"
means ``a pencil".) To make our consideration
clearer, we study this case twice: by usual
methods in Theorem \ref{t0.1} and by Arnold's
method of bifurcation diagrams in Theorem
\ref{t2.0}.

The theory of matrix pencils is the theory of
pairs of matrices of the same size up to
equivalence. Two matrix pairs $(A_1,A_2)$ and
$(B_1,B_2)$ are {\it equivalent} if there exist
two nonsingular matrices $R$ and $S$ such that
\begin{equation}\label{0}
  B_1=RA_1S,\qquad B_2=RA_2S.
\end{equation}
As was proved by Kronecker (see \cite{gan}),
every pair of complex $m\times n$ matrices is
equivalent to a direct sum, determined uniquely
up to permutation of summands, of pairs of the
form
\begin{equation}       \label{1.1}
(I_r,J_r(\lambda)),\ (J_r(0),I_r),\ (F_r,K_r),\
(F_r^T,K_r^T),
\end{equation}
where $J_r(\lambda)$ is the Jordan cell with
units over the diagonal and
\begin{equation*}       
F_r=\begin{bmatrix}
         1&&0\\
         0&\ddots&\\
         &\ddots&1\\
         0&&0
         \end{bmatrix},\qquad
K_r=\begin{bmatrix}
         0&&0\\
         1&\ddots&\\
         &\ddots&0\\
         0&&1
         \end{bmatrix}
\end{equation*}
are matrices of size  $r\times (r-1)$, $r\ge 1$.

Every square matrix is transformed to a
diagonalizable matrix by a small jiggling. This
result is extended to matrix pencils in the next
theorem (the same extension but in terms of
bifurcation diagrams is given in Theorem
\ref{t2.0}).

\begin{theorem} \label{t0.1}
Every pair $(A,B)$ of complex $m\times n$
matrices is transformed by an arbitrarily small
jiggling to a pair that is equivalent to the
following  pair:
\begin{itemize}

\item[(i)]
$(I_n, \diag(\alpha_1,\dots,\alpha_n))$ if $m=n$,

\item[(ii)]
$([I\ 0],[0\ I])$ if $m<n$, and

\item[(iii)]
$([I\ 0]^T,[0\ I]^T)$ if $m>n$.
\end{itemize}
\end{theorem}

\begin{proof}
Let $m=n$. We make $A$ nonsingular by an
arbitrarily small perturbation, reduce it to $I$,
then make $B$ diagonalizable by an arbitrarily
small perturbation.

Let $m\ne n$, suppose $m<n$ (otherwise, consider
the pair $(A^T,B^T)$). We make the rows of $A$
linearly independent by an arbitrarily small
perturbation, then reduce it to the form $A=[I\
0]$ and partition $B=[B_1\ B_2]$ conformal with
$A$.

If the number of rows of $B_2$ is no more than
the number of columns, we make its rows linearly
independent by an arbitrarily small perturbation,
reduce $B_2$ to the form $[0\ I]$ and make
$B_1=0$ by column transformations. The pair
$(A,B)$ takes the form (ii).

If the number of rows of $B_2$ is more than the
number of columns, we make the columns of $B_2$
linearly independent by an arbitrarily small
perturbation, reduce $B_2$ to the form $[0\
I]^T$, then reduce $B$ to the form
$$
B=\begin{bmatrix}
  B_{11}&0\\0&I
\end{bmatrix}
$$
by additions of columns of $B_2$. The
corresponding horizontal division of $A$ crosses
the block $I$, so we make a vertical division and
obtain
\begin{equation}\label{0.1}
A=\begin{bmatrix}
  I&0&0\\0&I&0
\end{bmatrix},\qquad
B=\begin{bmatrix}
  C_1&C_2&0\\0&0&I
\end{bmatrix}.
\end{equation}

Applying the same transformations to the fragment
$[I\ 0]$, $[C_1\ C_2]$, we reduce $(A,\ B)$
respectively to the form (ii) (making $C_2=[0\
I]$ and $C_1=0$) or to the form
$$
A=\begin{bmatrix}
  I&0&0&0\\0&I&0&0\\0&0&I&0
\end{bmatrix},\qquad
B=\begin{bmatrix}
  D_1&D_2&0&0\\0&0&I&0\\0&0&0&I
\end{bmatrix}
$$
(these transformations spoil the reduced part of
the pair \eqref{0.1}, but it is recovered by
obvious transformations). We repeat this
reduction until obtain $(A, B)$ of the form (ii).
\end{proof}

Note that miniversal deformations and bifurcation
diagrams for real matrices up to similarity were
given by Galin \cite{gal1}; for certain classes
of operators in metric spaces in \cite{ben},
\cite{gal2}, \cite{pat1}, \cite{pat2}, and
\cite{pat3}.

\section{Generic families of matrix pencils} \label{s1}

We study  families of pairs of complex $m\times
n$ matrices ${\cal
A}(\vec\alpha)=(A_1(\vec\alpha),A_2(\vec\alpha))$,
$\vec\alpha=(\alpha_1,\dots,\alpha_k)$,
holomorphic at $\vec 0$. The entries of
$A_1(\vec\alpha)$ and $A_2(\vec\alpha)$ are power
series of complex parameters $\alpha_1,\dots,
\alpha_k$ that are convergent in a neighborhood
of $\vec 0$. (The germ of a family ${\cal
A}(\vec\alpha)$ at $\vec 0$ is called a {\it
deformation} of the pair ${\cal A}(\vec 0)$, see
\cite{arn1}--\cite{arn2}.)

Two families ${\cal A}(\vec\alpha)$ and ${\cal
B}(\vec\alpha)$ are called {\it equivalent} if
there exist matrices $R(\vec\alpha)$ and
$S(\vec\alpha)$ holomorphic at $\vec 0$ such that
$R(\vec 0)=I$, $S(\vec 0)=I$, and
$$B_1(\vec\alpha)= R(\vec\alpha)
A_1(\vec\alpha)S(\vec\alpha), \qquad
B_2(\vec\alpha)= R(\vec\alpha)
A_2(\vec\alpha)S(\vec\alpha)$$ in a neighborhood
of $\vec 0$ (compare with \eqref{0}). A family
${\cal A}(\alpha_1,\dots,\alpha_k)$ is called
{\it versal} if every family ${\cal
B}(\beta_1,\dots,\beta_l)$ with ${\cal B}(\vec
0)= {\cal A}(\vec 0)$ is equivalent to a family
${\cal A}(\varphi_1(\vec\beta),\dots,
\varphi_k(\vec\beta))$, where
$\varphi_i(\vec\beta)$, $\varphi_i(\vec 0)=\vec
0$, are power series convergent in a neighborhood
of $\vec 0$. A versal family with the minimum
possible number $k$ of parameters is said to be
{\it miniversal}.

For every pair of $m\times n$ matrices
$(A_1,A_2)$, a miniversal family ${\cal
A}(\vec\alpha)$ with ${\cal A}(\vec 0)=(A_1,A_2)$
was obtained in \cite{kag} and simplified in
\cite{gar_ser}. We now  recall the result of
\cite{gar_ser}. It suffices to construct a
miniversal family for a Kronecker canonical pair
\begin{equation}       \label{1.3}
(A_1,A_2)=\bigoplus_{i=1}^l(F_{p_i},
K_{p_i})\oplus (I,C) \oplus (D,I)\oplus
\bigoplus_{i=1}^r(F_{q_i}^T, K_{q_i}^T)
\end{equation}
(see \eqref{1.1}), where $p_1\le\dots\le p_l$,
$q_1\ge\dots\ge q_r$,
\begin{equation} \label{1.4}
C=\bigoplus_{i=1}^t\Phi_i({\lambda_i})\quad
(\lambda_i\ne\lambda_j\ {\rm if}\ i\ne j),\qquad
D=\Phi_0(0),
\end{equation}
and $\Phi_0(0),\Phi_1({\lambda_1}),\dots,
\Phi_t({\lambda_t})$ have the form
$$
 \Phi_i({\lambda_i})={\rm
diag}(J_{s_{i1}}({\lambda_i}),\,
J_{s_{i2}}({\lambda_i}),\dots),\qquad s_{i1}\ge
s_{i2}\ge\cdots,\quad \lambda_0:=0.
$$

Denote by $0^{\uparrow}$ (resp.,
$0^{\downarrow},\ 0^{\leftarrow},\,
0^{\rightarrow}$) a matrix, in which all entries
are zero except for the entries of the first row
(respectively, the last row, the first column,
the last column) that are independent parameters;
and denote by $Z$ the $p\times q$ matrix, in
which the first $\max\{q-p,0\}$ entries of the
first row are independent parameters and the
other entries are zeros:
\begin{equation*}
0^{\uparrow}=\begin{bmatrix}
         *&\cdots&*\\
         0&\cdots&0\\
         \hdotsfor{3}\\
          0&\cdots&0
         \end{bmatrix},\qquad
Z=\left[\begin{tabular}{cccccc}
         $*$ & $\cdots$ & $*$ & 0 & $\cdots$ & 0\\
             &          &     &   & $\ddots$ &  \\
             & \LARGE 0 &     & 0 & $\cdots$ & 0
         \end{tabular}\right],
\end{equation*}
where the stars denote independent parameters.
Let
\begin{equation}\label{1.6}
{\cal H}=[H_{ij}]
\end{equation}
be a block matrix, whose $p_i\times q_j$ blocks
$H_{ij}$ are of the form $H_{ij}=0^{\leftarrow}$
if $p_i\le q_j$, and $H_{ij}=0^{\downarrow}$ if
$p_i>q_j$.

\begin{theorem}[see \cite{gar_ser}]   \label{t1.1}
One of miniversal families with the pair
\eqref{1.3} at $\vec 0$ is $\cal M(\vec\alpha)=$

{\rm $\!\!\!\!\!\!\!\!\left(
\special{em:linewidth 0.4pt} \unitlength 0.60mm
\linethickness{0.4pt}
\begin{picture}(212.00,60.00)(29,50)
\emline{130.00}{60.00}{1}{30.00}{60.00}{2}
\emline{90.00}{100.00}{3}{90.00}{0.00}{4}
\emline{80.00}{100.00}{5}{80.00}{40.00}{6}
\emline{80.00}{40.00}{7}{130.00}{40.00}{8}
\emline{130.00}{50.00}{9}{70.00}{50.00}{10}
\emline{70.00}{50.00}{11}{70.00}{100.00}{12}
\put(35.00,95.00){\makebox(0,0)[cc]{$F_{p_1}$}}
\put(45.00,85.00){\makebox(0,0)[cc]{$F_{p_2}$}}
\put(55.00,75.00){\makebox(0,0)[cc]{$\ddots$}}
\put(65.00,65.00){\makebox(0,0)[cc]{$F_{p_l}$}}
\put(75.00,80.00){\makebox(0,0)[cc]{\Large 0}}
\put(85.00,55.00){\makebox(0,0)[cc]{\Large 0}}
\put(110.00,55.00){\makebox(0,0)[cc]{\Large 0}}
\put(60.00,30.00){\makebox(0,0)[cc]{\Large 0}}
\put(75.00,55.00){\makebox(0,0)[cc]{$I$}}
\put(85.00,45.00){\makebox(0,0)[cc]{$\tilde D$}}
\put(60.00,90.00){\makebox(0,0)[cc]{\Large 0}}
\put(40.00,70.00){\makebox(0,0)[cc]{\Large 0}}
\put(95.00,35.00){\makebox(0,0)[cc]{$F_{q_1}^T$}}
\put(105.00,25.00){\makebox(0,0)[cc]{$F_{q_2}^T$}}
\put(115.00,15.00){\makebox(0,0)[cc]{$\ddots$}}
\put(125.00,5.00){\makebox(0,0)[cc]{$F_{q_r}^T$}}
\put(120.00,30.00){\makebox(0,0)[cc]{\Large 0}}
\put(100.00,10.00){\makebox(0,0)[cc]{\Large 0}}
\put(85.00,95.00){\makebox(0,0)[cc]{$0^{\downarrow}$}}
\put(85.00,85.00){\makebox(0,0)[cc]{$0^{\downarrow}$}}
\put(85.00,75.00){\makebox(0,0)[cc]{$\vdots$}}
\put(85.00,65.00){\makebox(0,0)[cc]{$0^{\downarrow}$}}
\put(95.00,80.00){\makebox(0,0)[cc]{$0^{\rightarrow}$}}
\put(105.00,80.00){\makebox(0,0)[cc]{$0^{\rightarrow}$}}
\put(115.00,80.00){\makebox(0,0)[cc]{$\cdots$}}
\put(125.00,80.00){\makebox(0,0)[cc]{$0^{\rightarrow}$}}
\put(125.00,45.00){\makebox(0,0)[cc]{$0^{\rightarrow}$}}
\put(115.00,45.00){\makebox(0,0)[cc]{$\cdots$}}
\put(105.00,45.00){\makebox(0,0)[cc]{$0^{\rightarrow}$}}
\put(95.00,45.00){\makebox(0,0)[cc]{$0^{\rightarrow}$}}
\put(135,50){\makebox(0,0)[cc]{,}}
\put(145.00,95.00){\makebox(0,0)[cc]{$K_{p_1}$}}
\put(155.00,85.00){\makebox(0,0)[cc]{$K_{p_2}$}}
\put(165.00,75.00){\makebox(0,0)[cc]{$\ddots$}}
\put(175.00,65.00){\makebox(0,0)[cc]{$K_{p_l}$}}
\put(195.00,55.00){\makebox(0,0)[cc]{\Large 0}}
\put(170.00,30.00){\makebox(0,0)[cc]{\Large 0}}
\put(185.00,55.00){\makebox(0,0)[cc]{$\tilde C$}}
\put(195.00,45.00){\makebox(0,0)[cc]{$I$}}
\put(150.00,70.00){\makebox(0,0)[cc]{\Large 0}}
\put(205.00,35.00){\makebox(0,0)[cc]{$K_{q_1}^T$}}
\put(215.00,25.00){\makebox(0,0)[cc]{$K_{q_2}^T$}}
\put(225.00,15.00){\makebox(0,0)[cc]{$\ddots$}}
\put(235.00,5.00){\makebox(0,0)[cc]{$K_{q_r}^T$}}
\put(210.00,10.00){\makebox(0,0)[cc]{\Large 0}}
\put(155.00,95.00){\makebox(0,0)[cc]{$Z$}}
\put(165.00,95.00){\makebox(0,0)[cc]{$\cdots$}}
\put(175.00,95.00){\makebox(0,0)[cc]{$Z$}}
\put(175.00,85.00){\makebox(0,0)[cc]{$\vdots$}}
\put(175.00,75.00){\makebox(0,0)[cc]{$Z$}}
\put(165.00,85.00){\makebox(0,0)[cc]{$\ddots$}}
\put(215.00,35.00){\makebox(0,0)[cc]{$Z^T$}}
\put(225.00,35.00){\makebox(0,0)[cc]{$\cdots$}}
\put(235.00,35.00){\makebox(0,0)[cc]{$Z^T$}}
\put(235.00,25.00){\makebox(0,0)[cc]{$\vdots$}}
\put(235.00,15.00){\makebox(0,0)[cc]{$Z^T$}}
\put(225.00,25.00){\makebox(0,0)[cc]{$\ddots$}}
\put(185.00,95.00){\makebox(0,0)[cc]{$0^{\uparrow}$}}
\put(185.00,85.00){\makebox(0,0)[cc]{$0^{\uparrow}$}}
\put(185.00,75.00){\makebox(0,0)[cc]{$\vdots$}}
\put(185.00,65.00){\makebox(0,0)[cc]{$0^{\uparrow}$}}
\put(220.00,95.00){\makebox(0,0)[cc]{$0^{\uparrow}$}}
\put(220.00,85.00){\makebox(0,0)[cc]{$0^{\uparrow}$}}
\put(220.00,75.00){\makebox(0,0)[cc]{$\vdots$}}
\put(220.00,65.00){\makebox(0,0)[cc]{$0^{\uparrow}$}}
\put(205.00,55.00){\makebox(0,0)[cc]{$0^{\leftarrow}$}}
\put(215.00,55.00){\makebox(0,0)[cc]{$0^{\leftarrow}$}}
\put(225.00,55.00){\makebox(0,0)[cc]{$\cdots$}}
\put(235.00,55.00){\makebox(0,0)[cc]{$0^{\leftarrow}$}}
\put(195.00,80.00){\makebox(0,0)[cc]{\Large 0}}
\put(220.00,45.00){\makebox(0,0)[cc]{\Large 0}}
\emline{145.00}{0.00}{13}{140.00}{0.00}{14}
\emline{140.00}{0.00}{15}{140.00}{100.00}{16}
\emline{140.00}{100.00}{17}{145.00}{100.00}{18}
\emline{235.00}{100.00}{19}{240.00}{100.00}{20}
\emline{240.00}{100.00}{21}{240.00}{0.00}{22}
\emline{240.00}{0.00}{23}{235.00}{0.00}{24}
\emline{240.00}{60.00}{25}{140.00}{60.00}{26}
\emline{200.00}{100.00}{27}{200.00}{0.00}{28}
\emline{240.00}{40.00}{29}{190.00}{40.00}{30}
\emline{190.00}{40.00}{31}{190.00}{100.00}{32}
\emline{180.00}{100.00}{33}{180.00}{50.00}{34}
\emline{180.00}{50.00}{35}{240.00}{50.00}{36}
\emline{125.00}{0.00}{37}{130.00}{0.00}{38}
\emline{130.00}{0.00}{39}{130.00}{100.00}{40}
\emline{130.00}{100.00}{41}{125.00}{100.00}{42}
\emline{35.00}{100.00}{43}{30.00}{100.00}{44}
\emline{30.00}{100.00}{45}{30.00}{0.00}{46}
\emline{30.00}{0.00}{47}{35.00}{0.00}{48}
\end{picture}
\right),$}\\[4mm]

\noindent where $$\tilde C=
\bigoplus_{i=1}^t(\Phi_i({\lambda_i})+{\cal
H}_i),\qquad \tilde D=\Phi_0(0)+{\cal H}_0$$ (see
\eqref{1.4}), and ${\cal H}_i$ is of the form
\eqref{1.6}.
\end{theorem}

Note that $\tilde C$ and $\tilde D$ are
miniversal deformations of $C$ and $D$ under
similarity, which were given by Arnold
\cite{arn1}--\cite{arn2}.

Extending Arnold's notation from
\cite{arn1}--\cite{arn2}, we will denote the
pairs \eqref{1.1} by the symbols
\begin{equation*}
\lambda^r:=(I_r,J_r(\lambda)),\
\infty^r:=(J_r(0),I_r),\
\vartriangle^r:=(F_r,K_r),\
\triangledown^r\!:=(F_r^T,K_r^T)
\end{equation*}
and a Kronecker canonical pair of matrices by a
sequence of these symbols. The complex number
$\lambda$ and the symbol $\infty$ will be called
the {\it eigenvalues} of $(I_r,J_r(\lambda))$ and
$(J_r(0),I_r)$; the eigenvalues will be denoted
by small Greek letters. We will say that two
pairs of $m\times n$ matrices have the same {\it
Kronecker type} if their Kronecker canonical
forms differ only by the sets of distinct
eigenvalues; the Kronecker type will be given by
an unordered sequence of symbols
$\vartriangle^r,\ \triangledown^r,\ \lambda^r$
($r\in {\mathbb N},\ \lambda\in {\mathbb
C}\cup\infty$, the set of $\lambda$'s in the
sequence is determined up to bijections in
${\mathbb C}\cup\infty$). In particular,
$(I_r,J_r(5))$ and $(J_r(0),I_r)$ have the same
Kronecker type $\lambda^r$; $(I_2,J_2(5))$,
$(I_2, 5I_2)$, and $(I_2,\diag(5,6))$ have
distinct Kronecker types $\lambda^2$,
$\lambda\lambda$, and $\lambda\mu$.

In the next section, we will study the set of
Kronecker types of matrix pairs that form a
miniversal family $\cal M(\vec\alpha)$ from
Theorem \ref{t1.1} in a neighborhood of $\vec 0$.
In this section, we remove parameters that have
no effect on the Kronecker type. Let us denote by
\begin{equation}\label{a}
{\cal M}'(\vec\beta)
\end{equation}
the family that is obtained from $\cal
M(\vec\alpha)$ by replacement of its blocks
$\tilde C$ and $\tilde D$ with
$${\tilde C}'=
\bigoplus_i(\Phi_i({\lambda_i})+{\cal H}'_i)\quad
\text{and}\quad {\tilde D}'=\Phi_0(0)+{\cal
H}'_0,
$$
where ${\cal H}'_i$ is obtained from
${\cal H}_i$ by replacement of its
upper-left-hand entry with 0 (the number of
parameters decreases by the number of distinct
eigenvalues). Furthermore, denote by
\begin{equation}\label{b}
 {\cal M}''(\vec\beta, \vec\gamma)
\end{equation}
the family that is obtained from ${\cal
M}'(\vec\beta)$ by replacement of ${\tilde C}'$
with $${\tilde C}''=
\bigoplus_i(\Phi_i({\lambda_i})+{\cal H}'_i+
\gamma_iI)=\bigoplus_i(\Phi_i({\lambda_i+\gamma_i})
+{\cal H}'_i)$$ and ${\tilde D}'$ with $${\tilde
D}''=\Phi_0(0)+{\cal H}'_0+
\gamma_0I=\Phi_0(\gamma_0)+{\cal H}'_0.$$

The families ${\cal M}''(\vec\beta, \vec\gamma)$
and $\cal M(\vec\alpha)$ have the same number of
parameters. Moreover, ${\cal M}''(\vec\beta,
\vec\gamma)$ is a miniversal family too; this
fact is proved in the same way as the
miniversality of $\cal M(\vec\alpha)$ (Theorem
2.1 from \cite{gar_ser}) and is based on the
following criterion.

{\it A criterion of miniversality} (see
\cite{arn1}--\cite{arn2} and \cite{gar_ser}): A
family
\begin{equation}\label{1.10}
{\cal A}(\alpha_1,\dots,\alpha_t)= {\cal A}_0
+\sum_{i=1}^t\alpha_i{\cal A}_i+...\,,\quad {\cal
A}_i=(A_i, B_i)\in {\mathbb C}^{(m\times
n,m\times n)},
\end{equation}
(where ${\mathbb C}^{(m\times n,m\times n)}$ is
the vector space of pairs of complex $m\times n$
matrices, and the points after + denote the terms
of order more than 1) is miniversal if and only
if
$${\mathbb C}^{(m\times
n,m\times n)}={\cal P}_{\cal A} \oplus {\cal
T}_{\cal A},$$ where
\begin{equation}\label{1.11}
{\cal P}_{\cal A}=\{\alpha_1{\cal A}_1+\dots +
\alpha_t{\cal A}_t\,|\,\alpha_i\in\mathbb C\}
\end{equation}
is the vector space spanned by ${\cal A}_1,\dots,
{\cal A}_t$, and
\begin{equation*}\label{1.12}
{\cal T}_{\cal A}=\{(RA_0-A_0S,RB_0-B_0S)\,|\,
R\in {\mathbb C}^{m\times m},\ S\in {\mathbb
C}^{n\times n}\}
\end{equation*}
is the tangent space at the point ${\cal A}_0$ to
the equivalence class
$$\{(RA_0S,RB_0S)\,|\, R\in
\gl_m({\mathbb C}),\ S\in \gl_n({\mathbb C})\}$$
of the pair ${\cal A}_0$.

Similar to \cite{arn1}--\cite{arn2}, we say that
a family \eqref{1.10} is {\it transversal to the
stratification into Kronecker types} if
\begin{equation}\label{1.13}
{\mathbb C}^{(m\times n,m\times n)}={\cal
P}_{\cal A} + {\cal Q}_{\cal A},
\end{equation}
where ${\cal P}_{\cal A}$ is defined by
\eqref{1.11} and ${\cal Q}_{\cal A}$ is the
tangent space at the point ${\cal A}_0$ to the
class of all pairs of matrices having the same
Kronecker type as ${\cal A}_0$. (Two subspaces
are {\it transversal} if their sum is the entire
space.)

\begin{theorem} \label{t1.2}
(i) In the space of families of pairs of $m\times
n$ matrices, the families transversal to the
stratification into Kronecker types constitute an
everywhere dense set.

(ii) For every miniversal family ${\cal
M}(\vec\alpha)$ from Theorem \ref{t1.1}, the
family ${\cal M}'(\vec\beta)$ (see \eqref{a}) is
transversal to the stratification into Kronecker
types. Moreover, the corresponding sum
\eqref{1.13} is direct:
\begin{equation}\label{1.14}
{\mathbb C}^{(m\times n,m\times n)}= {\cal
P}_{{\cal M}'}\oplus {\cal Q}_{{\cal M}'}.
\end{equation}
\end{theorem}

\begin{proof}
The statement (i) follows from  the theorem of
\cite[\S{\,}30E]{arn2}. Let us prove the
statement (ii). The family \eqref{b}  has the
form
$${\cal
M}''(\vec\beta, \vec\gamma)={\cal M}'(\vec\beta)+
\gamma_1{\cal M}''_1+ \dots+\gamma_s{\cal
M}''_s$$ and has the same Kronecker type as
${\cal M}'(\vec\beta)$ for a small $\vec\gamma$,
therefore, $\gamma_1{\cal M}''_1+
\dots+\gamma_s{\cal M}''_s\in {\cal Q}_{{\cal
M}'}$.  Since ${\cal M}''(\vec\beta, \vec\gamma)$
is a miniversal family,
\begin{align*}
 {\mathbb C}^{(m\times n,m\times n)}&=
 {\cal P}_{{\cal M}''} \oplus {\cal T}_{{\cal
 M}''}\\
 &={\cal P}_{{\cal M}'} \oplus \{\gamma_1{\cal
M}''_1+ \dots+\gamma_s{\cal M}''_s\,|\,
\gamma_i\in \mathbb C\}\oplus {\cal T}_{{\cal
M}''}\\
 &={\cal P}_{{\cal M}'}\oplus {\cal
Q}_{{\cal M}'},
\end{align*}
this proves \eqref{1.14}.
\end{proof}

Similar to \cite{arn1}--\cite{arn2}, a family
will be called a {\it generic family} (or a
family in general position) if it is transversal
to the stratification into Kronecker types.

\begin{corollary} \label{col}
A nongeneric family can be transformed into a
generic family by an arbitrarily small
perturbation of the family. Since the sum
\eqref{1.14} is a direct sum, the families ${\cal
M}'(\vec\beta)$ have the most complicated
Kronecker structure among the generic families:
if an arbitrary family ${\cal A}(\vec\delta)$ has
the same number of parameters as ${\cal
M}'(\vec\beta)$ but contains matrices with more
complicated Kronecker structure, then ${\cal
A}(\vec\delta)$ is not a generic family.
\end{corollary}

\section{Bifurcation diagrams} \label{s2}

In this section, we construct bifurcation
diagrams for generic zero-, one- and
two-parameter families of pairs of matrices.

Let ${\cal A}(\vec\alpha)$ be a family of pairs
of ${m\times n}$ matrices; that is, a holomorphic
mapping $${\cal A}:\Lambda\to {\mathbb
C}^{(m\times n,m\times n)},$$ where
$\Lambda\subset{\mathbb C}^k$ is a neighborhood
of $\vec 0$. A {\it bifurcation diagram} of this
family is, by definition, a partition of the
parameter domain $\Lambda$ according to Kronecker
types of pairs. To construct it we assign to each
$\vec\alpha\in\Lambda$ the Kronecker type of
${\cal A}(\vec\alpha)$ and then join all points
with the same Kronecker type. We narrow down the
neighborhood $\Lambda$ when it simplifies the
structure of the bifurcation diagram. The
bifurcation diagram of a generic family reflects
the possible Kronecker structure of pairs in the
family. By \cite[\S{\,}30E]{arn2}, if in the
study of a phenomenon we obtain another
bifurcation diagram then  in the idealization of
the phenomenon something essential was missed, or
there were some special reasons for an additional
complexity of the structure, or the family is not
generic.

If all pairs ${\cal A}(\vec\alpha)$,
$\vec\alpha\in\Lambda$, have the same Kronecker
type $\bar t$, we will give the bifurcation
diagram by the sequence $\bar t$. But usually
matrix pairs of a generic family have distinct
types at $\vec 0$ and outside of $\vec 0$; in
this case we will give the bifurcation diagram by
the pair ${\bar t}_0/{\bar t}_1$, where ${\bar
t}_0$ is the Kronecker type of ${\cal A}(\vec 0)$
and ${\bar t} _1$ is the set of Kronecker types
of ${\cal A}(\vec\alpha)$, $\vec 0\ne
\vec\alpha\in\Lambda$.

\begin{theorem} \label{t2.0}
Generic zero-parameter families of pairs of
$m\times n$ matrices have
 the bifurcation diagrams
\begin{equation}\label{e1}
\left.\begin{tabular}{lc}
$\vartriangle^r\!\!\ldots\vartriangle^r
\vartriangle^{r+1}\!\!\ldots\vartriangle^{r+1}\qquad\qquad$
& $\triangledown^r\!\ldots\triangledown^r
\triangledown^{r+1}\!\ldots\triangledown^{r+1}$\\
$\lambda\mu\ldots\tau$ &
\end{tabular}\right\}
\end{equation}
(the parts
$\vartriangle^{r+1}\!\!\ldots\vartriangle^{r+1}$
and
$\triangledown^{r+1}\!\ldots\triangledown^{r+1}$
can be absent).
\end{theorem}

\begin{proof}
Let $\cal A$ be a generic zero-parameter family.
By Corollary \ref{col}, we may suppose that
${\cal A}={\cal M}'(\vec\beta)$, where ${\cal
M}(\vec\alpha)$ is a family from Theorem
\ref{t1.1}. Selecting all ${\cal M}(\vec\alpha)$
for which ${\cal M}'(\vec\beta)$ is a
zero-parameter family, we obtain the list
\eqref{e1} of all possible Kronecker types for
$\cal A$. (Note that the pairs (i)--(iii) from
Theorem \ref{t0.1} have the types \eqref{e1};
this gives another proof of Theorem \ref{t0.1}.)
\end{proof}

\begin{theorem} \label{t2.1}
Generic one-parameter families of pairs of
$m\times n$ matrices have the bifurcation
diagrams \eqref{e1}, in which case the families
behave as generic zero-parameter families, or
have the following bifurcation diagrams:
\begin{equation}\label{e2}
\left.\begin{tabular}{ll}
$\vartriangle^r\vartriangle^{r+2}\!/
\!\vartriangle^{r+1}\vartriangle^{r+1}$&
$\triangledown^r\triangledown^{r+2}\,/\,
\triangledown^{r+1}\triangledown^{r+1}$\\
$\vartriangle^r\!\!\lambda\,/\!\vartriangle^{r+1}$&
$\triangledown^r\lambda\,/\,\triangledown^{r+1}$\\
$\lambda^2_1\lambda_2\ldots\lambda_t\,/\,
\mu_1\mu_2\ldots\mu_{t+1}\qquad$&
\end{tabular}\right\}
\end{equation}
If in a one-parameter family there are pairs with
a more complicated Kronecker structure, then we
can remove them by an arbitrarily small
perturbation of the family.
\end{theorem}

\begin{proof}
Let ${\cal A}(\beta)$ be a generic one-parameter
family. If \eqref{1.13} is not a direct sum, then
the parameter $\beta$ does not affect on the
Kronecker type and all pairs of the family have
the same type. Applying Theorem \ref{t2.0}, we
get the list \eqref{e1} of admissible bifurcation
diagrams for the family.

Suppose \eqref{1.13} is a direct sum. By
Corollary \ref{col}, we may take ${\cal
A}(\beta)={\cal M}'(\beta)$, where ${\cal
M}(\vec\alpha)$ is a family from Theorem
\ref{t1.1}. Selecting all ${\cal M}(\vec\alpha)$
for which ${\cal M}'(\vec{\beta})$ is a
one-parameter family, we obtain that ${\cal
A}(0)={\cal M}'(0)$ is one of the pairs:
\begin{equation*}\label{2.1}
\vartriangle^r\vartriangle^{r+2},\quad
\triangledown^r\triangledown^{r+2},\quad
\vartriangle^r\!\!\lambda,\quad
\triangledown^r\lambda,\quad
\lambda^2_1\lambda_2\ldots\lambda_t.
\end{equation*}

1) Let ${\cal A}(0)=
\vartriangle^r\vartriangle^{r+2}= (F_r,K_r)\oplus
(F_{r+2},K_{r+2})$ (see \eqref{1.1}). Then ${\cal
A}(\beta)$ is the pair of matrices
$$
\unitlength 0.70mm \linethickness{0.4pt}
\begin{picture}(155.00,85.00)
\put(10.00,0.00){\framebox(60.00,80.00)[cc]{}}
\put(70.00,40.00){\line(-1,0){60.00}}
\put(10.00,40.00){\line(0,0){0.00}}
\put(10.00,40.00){\line(0,0){0.00}}
\put(40.00,0.00){\line(0,1){80.00}}
\put(5.00,75.00){\makebox(3,0)[rc]{$\scriptstyle{1}$}}
\put(5.00,5.00){\makebox(3,0)[rc]{$\scriptstyle
\overline{r+2}$}}
\put(5.00,65.00){\makebox(3,0)[rc]{$\scriptstyle
2$}}
\put(5.00,45.00){\makebox(3,0)[rc]{$\scriptstyle
r$}}
\put(5.00,35.00){\makebox(3,0)[rc]{$\scriptstyle
\overline{1}$}}
\put(5.00,25.00){\makebox(3,0)[rc]{$\scriptstyle
\overline{2}$}}
\put(15.00,85.00){\makebox(0,-3)[cc]{$\scriptstyle
1$}}
\put(35.00,85.00){\makebox(0,-3)[cc]{$\scriptstyle{r-1}$}}
\put(45.00,85.00){\makebox(0,-3)[cc]{$\scriptstyle
\overline{1}$}}
\put(65.00,85.00){\makebox(0,-3)[cc]{$\scriptstyle
\overline{r+1}$ }}
\put(15.00,75.00){\makebox(0,0)[cc]{1}}
\put(15.00,65.00){\makebox(0,0)[cc]{0}}
\put(35.00,45.00){\makebox(0,0)[cc]{0}}
\put(35.00,55.00){\makebox(0,0)[cc]{1}}
\put(25.00,65.00){\makebox(0,0)[cc]{$\ddots$}}
\put(25.00,55.00){\makebox(0,0)[cc]{$\ddots$}}
\put(45.00,35.00){\makebox(0,0)[cc]{1}}
\put(45.00,25.00){\makebox(0,0)[cc]{0}}
\put(65.00,5.00){\makebox(0,0)[cc]{0}}
\put(65.00,15.00){\makebox(0,0)[cc]{1}}
\put(55.00,25.00){\makebox(0,0)[cc]{$\ddots$}}
\put(55.00,15.00){\makebox(0,0)[cc]{$\ddots$}}
\put(55.00,60.00){\makebox(0,0)[cc]{\LARGE 0 }}
\put(25.00,20.00){\makebox(0,0)[cc]{\LARGE 0}}
\put(5.00,15.00){\makebox(3,0)[rc]{$\vdots$}}
\put(5.00,55.00){\makebox(3,0)[rc]{$\vdots$}}
\put(25.00,85.00){\makebox(0,-3)[cc]{$\cdots$}}
\put(55.00,85.00){\makebox(0,-3)[cc]{$\cdots$}}
\put(95.00,0.00){\framebox(60.00,80.00)[cc]{}}
\put(155.00,40.00){\line(-1,0){60.00}}
\put(95.00,40.00){\line(0,0){0.00}}
\put(95.00,40.00){\line(0,0){0.00}}
\put(125.00,0.00){\line(0,1){80.00}}
\put(90.00,75.00){\makebox(3,0)[rc]{$\scriptstyle{1}$}}
\put(90.00,5.00){\makebox(3,0)[rc]{$\scriptstyle
\overline{r+2}$}}
\put(90.00,65.00){\makebox(3,0)[rc]{$\scriptstyle
2$}}
\put(90.00,45.00){\makebox(3,0)[rc]{$\scriptstyle
r$}}
\put(90.00,35.00){\makebox(3,0)[rc]{$\scriptstyle
\overline{1}$}}
\put(90.00,25.00){\makebox(3,0)[rc]{$\scriptstyle
\overline{2}$}}
\put(100.00,85.00){\makebox(0,-3)[cc]{$\scriptstyle
1$}}
\put(120.00,85.00){\makebox(0,-3)[cc]{$\scriptstyle{r-1}$}}
\put(130.00,85.00){\makebox(0,-3)[cc]{$\scriptstyle
\overline{1}$}}
\put(150.00,85.00){\makebox(0,-3)[cc]{$\scriptstyle
\overline{r+1}$ }}
\put(100.00,75.00){\makebox(0,0)[cc]{0}}
\put(100.00,65.00){\makebox(0,0)[cc]{1}}
\put(120.00,45.00){\makebox(0,0)[cc]{1}}
\put(120.00,55.00){\makebox(0,0)[cc]{0}}
\put(110.00,65.00){\makebox(0,0)[cc]{$\ddots$}}
\put(110.00,55.00){\makebox(0,0)[cc]{$\ddots$}}
\put(130.00,35.00){\makebox(0,0)[cc]{0}}
\put(130.00,25.00){\makebox(0,0)[cc]{1}}
\put(150.00,5.00){\makebox(0,0)[cc]{1}}
\put(150.00,15.00){\makebox(0,0)[cc]{0}}
\put(140.00,25.00){\makebox(0,0)[cc]{$\ddots$}}
\put(140.00,15.00){\makebox(0,0)[cc]{$\ddots$}}
\put(140.00,60.00){\makebox(0,0)[cc]{\LARGE 0 }}
\put(110.00,20.00){\makebox(0,0)[cc]{\LARGE 0}}
\put(90.00,15.00){\makebox(3,0)[rc]{$\vdots$}}
\put(90.00,55.00){\makebox(3,0)[rc]{$\vdots$}}
\put(110.00,85.00){\makebox(0,-3)[cc]{$\cdots$}}
\put(140.00,85.00){\makebox(0,-3)[cc]{$\cdots$}}
\put(130.00,75.00){\makebox(0,0)[cc]{$\beta$}}
\end{picture}
$$
(we always number the rows and columns of the
second summand by overbarred natural numbers).
Rearranging the rows and columns in the order
$$(\overline{1},\overline{2}\,|\,1,\overline{3}\,|\,2,
\overline{4}\,|\,\ldots\,|\,r-1,\overline{r+1}\,|\,r,
\overline{r+2})$$ and
$$(\overline{1},\overline{2}\,|\,1,\overline{3}\,|\,2,
\overline{4}\,|\,\ldots\,|\,r-2,\overline{r}\,|\,r-1,
\overline{r+1})$$ (the symbol $|$ denotes the
partition into strips), we obtain the following
pair of matrices:\\
$$
\unitlength 0.65mm \linethickness{0.4pt}
\begin{picture}(195.00,105.00)
\put(10.00,0.00){\framebox(80.00,100.00)[cc]{}}
\put(90.00,80.00){\line(-1,0){80.00}}
\put(10.00,60.00){\line(1,0){80.00}}
\put(90.00,40.00){\line(-1,0){80.00}}
\put(10.00,20.00){\line(1,0){80.00}}
\put(70.00,0.00){\line(0,1){100.00}}
\put(50.00,100.00){\line(0,-1){100.00}}
\put(30.00,0.00){\line(0,1){100.00}}
\put(15.00,105.00){\makebox(0,-3)[cc]{$\scriptstyle
\overline{1}$}}
\put(25.00,105.00){\makebox(0,-3)[cc]{$\scriptstyle
\overline{2}$}}
\put(35.00,105.00){\makebox(0,-3)[cc]{$\scriptstyle
1$}}
\put(45.00,105.00){\makebox(0,-3)[cc]{$\scriptstyle
\overline{3}$}}
\put(60.00,105.00){\makebox(0,-3)[cc]{$\cdots$}}
\put(75.00,105.00){\makebox(0,-3)[cc]{$\scriptstyle
r-1$}}
\put(85.00,105.00){\makebox(0,-2)[cc]{$\scriptstyle
\overline{r+1}$}}
\put(5.00,95.00){\makebox(3,0)[rc]{$\scriptstyle
\overline{1}$}}
\put(5.00,85.00){\makebox(3,0)[rc]{$\scriptstyle
\overline{2}$}}
\put(5.00,75.00){\makebox(3,0)[rc]{$\scriptstyle
1$}}
\put(5.00,65.00){\makebox(3,0)[rc]{$\scriptstyle
\overline{3}$}}
\put(5.00,51.00){\makebox(3,0)[rc]{$\vdots$}}
\put(5.00,35.00){\makebox(3,0)[rc]{$\scriptstyle
r-1$}}
\put(5.00,25.00){\makebox(3,0)[rc]{$\scriptstyle
\overline{r+1}$}}
\put(5.00,15.00){\makebox(3,0)[rc]{$\scriptstyle
r$}}
\put(5.00,5.00){\makebox(3,0)[rc]{$\scriptstyle
\overline{r+2}$}}
\put(115.00,0.00){\framebox(80.00,100.00)[cc]{}}
\put(195.00,80.00){\line(-1,0){80.00}}
\put(115.00,60.00){\line(1,0){80.00}}
\put(195.00,40.00){\line(-1,0){80.00}}
\put(115.00,20.00){\line(1,0){80.00}}
\put(175.00,0.00){\line(0,1){100.00}}
\put(155.00,100.00){\line(0,-1){100.00}}
\put(135.00,0.00){\line(0,1){100.00}}
\put(120.00,105.00){\makebox(0,-3)[cc]{$\scriptstyle
\overline{1}$}}
\put(130.00,105.00){\makebox(0,-3)[cc]{$\scriptstyle
\overline{2}$}}
\put(145.00,105.00){\makebox(0,-3)[cc]{$\cdots$}}
\put(160.00,105.00){\makebox(0,-3)[cc]{$\scriptstyle
r-2$}}
\put(170.00,105.00){\makebox(0,-2.5)[cc]{$\scriptstyle
\overline{r}$}}
\put(180.00,105.00){\makebox(0,-3)[cc]{$\scriptstyle
r-1$}}
\put(190.00,105.00){\makebox(0,-2)[cc]{$\scriptstyle
\overline{r+1}$}}
\put(110.00,95.00){\makebox(3,0)[rc]{$\scriptstyle
\overline{1}$}}
\put(110.00,85.00){\makebox(3,0)[rc]{$\scriptstyle
\overline{2}$}}
\put(110.00,75.00){\makebox(3,0)[rc]{$\scriptstyle
1$}}
\put(110.00,65.00){\makebox(3,0)[rc]{$\scriptstyle
\overline{3}$}}
\put(110.00,51.00){\makebox(3,0)[rc]{$\vdots$}}
\put(110.00,35.00){\makebox(3,0)[rc]{$\scriptstyle
r-1$}}
\put(110.00,25.00){\makebox(3,0)[rc]{$\scriptstyle
\overline{r+1}$}}
\put(110.00,15.00){\makebox(3,0)[rc]{$\scriptstyle
r$}}
\put(110.00,5.00){\makebox(3,0)[rc]{$\scriptstyle
\overline{r+2}$}}
\put(15.00,95.00){\makebox(0,0)[cc]{1}}
\put(25.00,85.00){\makebox(0,0)[cc]{1}}
\put(35.00,75.00){\makebox(0,0)[cc]{1}}
\put(45.00,65.00){\makebox(0,0)[cc]{1}}
\put(60.00,51.00){\makebox(0,0)[cc]{$\ddots$}}
\put(75.00,35.00){\makebox(0,0)[cc]{1}}
\put(85.00,25.00){\makebox(0,0)[cc]{1}}
\put(120.00,95.00){\makebox(0,0)[cc]{0}}
\put(130.00,95.00){\makebox(0,0)[cc]{0}}
\put(130.00,85.00){\makebox(0,0)[cc]{0}}
\put(120.00,85.00){\makebox(0,0)[cc]{1}}
\put(120.00,75.00){\makebox(0,0)[cc]{$\beta$}}
\put(130.00,75.00){\makebox(0,0)[cc]{0}}
\put(130.00,65.00){\makebox(0,0)[cc]{1}}
\put(120.00,65.00){\makebox(0,0)[cc]{0}}
\put(145.00,51.00){\makebox(0,0)[cc]{$\ddots$}}
\put(160.00,35.00){\makebox(0,0)[cc]{1}}
\put(170.00,25.00){\makebox(0,0)[cc]{1}}
\put(180.00,15.00){\makebox(0,0)[cc]{1}}
\put(190.00,5.00){\makebox(0,0)[cc]{1}}
\end{picture}
$$\\[-5 mm]

We may reduce this pair by simultaneous
elementary transformations. Let us prove that a
linear combination of rows of the (2,1) block of
the second matrix may be added to rows of its
(1,1) block without spoiling the other blocks of
the pair. Indeed, additions of rows of the second
horizontal strip to rows of the first horizontal
strip spoil the (1,2) block of the first matrix.
We recover it by additions of columns of the
first vertical strip to columns of the second
vertical strip, which spoil the (1,2) and (2,2)
blocks of the second matrix. We recover them by
additions of rows of the third horizontal strip
spoiling the (1,3) and (2,3) blocks of the first
matrix. These blocks are recovered by additions
of columns of the first and the second vertical
strips, and so on. On the last step, we recover
blocks of the last vertical strip of the second
matrix by additions of rows of the last
horizontal strip without spoiling the other
blocks since the last horizontal strip of the
first matrix is zero.

Let the parameter $\beta\ne 0$. Multiplying the
third rows by $\beta^{-1}$, we make unit the
(3,1) entry of the second matrix spoiling the
(3.3) entry of the first matrix. We recover it
multiplying the third columns by $\beta$ and so
on until obtain the initial pair with $\beta=1$.
The subtraction of the first row of the (2,1)
block of the second matrix from the second row of
its (1,1) block makes zero this block. Up to
simultaneous permutation of rows and columns, the
obtained pair has the form
$(F_{r+1},K_{r+1})\oplus (F_{r+1},K_{r+1})$.
Therefore, the bifurcation diagram of ${\cal
A}(\beta)$ is
$\vartriangle^r\vartriangle^{r+2}\!\!/
\!\vartriangle^{r+1}\vartriangle^{r+1}$.
\medskip

2) The case ${\cal A}(0)=
\triangledown^r\triangledown^{r+2}$ is considered
analogously.
\medskip

3) Let ${\cal A}(0)= \vartriangle^r\!\!\lambda$.
First, suppose $\lambda\ne \infty$. Up to
simultaneous permutation of rows and columns,
${\cal A}(\beta)$ is the pair $$
\begin{tabular}{|c|ccc|} \hline
 1&&&\\ \hline
 &1&&\\
 &0&$\ddots$&\\
 &&$\ddots$&1\\
 &&&0\\ \hline
\end{tabular}\qquad
\begin{tabular}{|c|ccc|} \hline
 $\lambda$ &&&\\ \hline
 $\beta$ &0&&\\
 &1&$\ddots$&\\
 &&$\ddots$&0\\
 &&&1\\ \hline
\end{tabular}
$$
Analogous to the case 1), we make $\beta=1$ in
the second matrix, then we make $\lambda=0$ in
this matrix by adding $\beta$. The obtained pair
is $(F_{r+1},K_{r+1})$.

Let now $\lambda= \infty$. Then ${\cal A}(\beta)$
is the pair
$$
\begin{tabular}{|ccc|c|} \hline
 1&&&\\
 0&$\ddots$&&\\
 &$\ddots$&1&\\
 &&0&$\beta$\\ \hline
 &&&0\\ \hline
\end{tabular}\qquad
\begin{tabular}{|ccc|c|} \hline
 0&&&\\
 1&$\ddots$&&\\
 &$\ddots$&0&\\
 &&1&\\ \hline
 &&&1\\ \hline
\end{tabular}
$$
Making $\beta=1$ gives the pair
$(F_{r+1},K_{r+1})$.

Therefore, the bifurcation diagram of ${\cal
A}(\beta)$ is
$\vartriangle^r\!\!\lambda\,/\!\vartriangle^{r+1}$.
\medskip

4) The case ${\cal A}(0)= \triangledown^r\lambda$
is considered analogously.
\medskip

5) Let ${\cal A}(0)=
\lambda^2_1\lambda_2\ldots\lambda_t$. First, we
suppose  $\lambda_1\ne\infty$. Then the pair
${\cal A}(\beta)$ has a direct summand
$$
\left(\begin{bmatrix}
  1&0 \\0&1
\end{bmatrix},\
\begin{bmatrix}
  \lambda_1&1 \\ \beta&\lambda_1
\end{bmatrix}\right).
$$
This pair reduces to $(1,\mu_1)\oplus(1,\mu_2)$
if $\beta\ne 0$. The case $\lambda_1=\infty$ is
considered similarly. The bifurcation diagram of
${\cal A}(\beta)$ is
$\lambda^2_1\lambda_2\ldots\lambda_t\,/\,
\mu_1\mu_2\ldots\mu_{t+1}$.
\end{proof}

\begin{theorem} \label{t2.2}
Generic two-parameter families of pairs of
$m\times n$ matrices have the bifurcation
diagrams  \eqref{e1} and \eqref{e2}, in which
case the families behave as generic
zero-parameter and one-parameter families, or
have the following bifurcation diagrams:
\begin{itemize}

\item[(i)]
$\vartriangle^1\!\triangledown^1\,/\,\lambda$,

\item[(ii)]
 $\vartriangle^r\vartriangle^{r+3}/
\!\vartriangle^{r+1}\vartriangle^{r+2}$,

\item[(iii)]
 $\vartriangle^r\vartriangle^r
\vartriangle^{r+2}\!/
\!\vartriangle^{r}\vartriangle^{r+1}
\vartriangle^{r+1}$,

\item[(iv)]
$\vartriangle^r\vartriangle^{r+2}
\vartriangle^{r+2}\!/
\!\vartriangle^{r+1}\vartriangle^{r+1}
\vartriangle^{r+2}$,

\item[(v)]
$\vartriangle^r\vartriangle^{r}\!\lambda\,/
\!\vartriangle^{r}\vartriangle^{r+1}$,

\item[(vi)]
$\vartriangle^r\vartriangle^{r+1}\!\!\lambda\,/\,
\{\vartriangle^r\vartriangle^{r+2}$ (this type
have the pairs with parameters on a smooth line
through $\vec 0$),
$\vartriangle^{r+1}\vartriangle^{r+1}$ (the pairs
with parameters outside the line)\},

\item[(vii)]
$\vartriangle^r\!\!\lambda\mu\,/\,
\{\vartriangle^{r+1}\!\!\nu$ (the pairs with
parameters on two smooth lines intersecting at
$\vec 0$), $\vartriangle^{r+2}$ (the pairs with
parameters outside the lines)\},

\item[(viii)]
$\lambda^3_1\lambda_2\ldots\lambda_t\,/\,
\{\mu^2_1\mu_2\ldots\mu_{t+1}$ (the pairs with
parameters on a line with a cusp at $\vec 0$),
$\mu_1\mu_2\ldots\mu_{t+2}$ (the pairs with
parameters outside the line)\},

\item[(ix)]
$\lambda^2_1\lambda^2_2\lambda_3
\ldots\lambda_t\,/\,
\{\mu^2_1\mu_2\ldots\mu_{t+1}$ (the pairs with
parameters on two smooth lines intersecting at
$\vec 0$), $\mu_1\mu_2\ldots\mu_{t+2}$ (the pairs
with parameters outside the lines)\},

\item[(x)]
the diagrams that are obtained from the diagrams
(ii)--(vii) by replacing all symbols
$\vartriangle$ by $\triangledown$.
\end{itemize}
If in a two-parameter family there are pairs with
a more complicated Kronecker structure, then we
can remove them by an arbitrarily small
perturbation of the family.
\end{theorem}

\begin{proof}
Let ${\cal A}(\beta,\gamma)$ be a generic
two-parameter family. If \eqref{1.13} is not a
direct sum, then the family behaves as a
one-parameter family, so its bifurcation diagram
is contained in the lists \eqref{e1} and
\eqref{e2}.

Let \eqref{1.13} be a direct sum, then we may
take ${\cal A}(\beta,\gamma)={\cal
M}'(\beta,\gamma)$, where ${\cal M}(\vec\alpha)$
is a family from Theorem \ref{t1.1}. Selecting
all ${\cal M}(\vec\alpha)$ for which ${\cal
M}'(\vec\beta)$ is a two-parameter family, we
obtain that ${\cal A}(0,0)={\cal M}'(0,0)$ is one
of the pairs:
\begin{equation}\label{2.1a}
\left.\begin{tabular}{lll}
$\vartriangle^1\!\triangledown^1$ &
$\vartriangle^r\vartriangle^{r+3}$ &
$\vartriangle^r\vartriangle^r \vartriangle^{r+2}$
\\
$\vartriangle^r\vartriangle^{r+2}
\vartriangle^{r+2}  \quad   \qquad$ &
$\vartriangle^r\vartriangle^{r}\!\lambda$ &
$\vartriangle^r\vartriangle^{r+1}\!\!\lambda$ \\
$\vartriangle^r\!\!\lambda\mu$ &
$\lambda^3_1\lambda_2\ldots\lambda_t \quad \qquad
$ &
$\lambda^2_1\lambda^2_2\lambda_3 \ldots\lambda_t$
\end{tabular}\right\}
\end{equation}
or it is obtained from them by turnover of all
$\vartriangle$ and $\triangledown$. We will
consider only the pairs \eqref{2.1a}, the others
are reduced to them by taking the transposed
matrices.\medskip

1) Let ${\cal A}(0,0)=
\vartriangle^1\!\!\triangledown^1=(0,0)$. Then
${\cal A}(\beta,\gamma)=(\beta,\gamma)$; we have
the bifurcation diagram (i). \medskip

 2) Let ${\cal
A}(0,0)= \vartriangle^r\vartriangle^{r+3}$. Let
the rows and columns of ${\cal A}(\beta,\gamma)$
be numbered by $1,2,\dots,r;\overline{1},
\overline{2},\dots, \overline{r+3}$ and
$1,2,\dots,r-1;\overline{1}, \overline{2},\dots,
\overline{r+2}$. Rearranging them in the order
$$(\overline{1},\overline{2},\overline{3}\,|\,1,
\overline{4}\,|\,2,\overline{5}\,|\,\ldots\,|\,
r-1,\overline{r+2}\,|\,r,\overline{r+3})$$ and
$$(\overline{1},\overline{2},\overline{3}
\,|\,1, \overline{4}\,|\,2,\overline{5}\,|\, \ldots\,|\,
r-2,\overline{r+1}\,|\,r-1,\overline{r+2}),$$ we
obtain the following pair of matrices:
$$
\begin{tabular}{|ccc|cc|c|cc|}
  \hline
  1&&&  &&    &         &\\
  &1&&  &&    &         &\\
  &&1&  &&    &         &\\
  \hline
  &&&  1&&    &         &\\
  &&&  &1&    &         &\\
  \hline
  &&&   && $\ddots$&    &\\
  \hline
  &&&   &&    &        1&\\
  &&&   &&    &        &1\\
  \hline
  &&&   &&    &         &\\
  &&&   &&    &         &\\
  \hline
\end{tabular}\qquad
\begin{tabular}{|ccc|c|cc|cc|}
  \hline
  0&0&0&  &    &&         &\\
  1&0&0&  &    &&         &\\
  0&1&0&  &    &&         &\\
  \hline
  $\beta$&$\gamma$&0& &    &&    &\\
  0&0&1&  &    &&         &\\
  \hline
  &&&   $\ddots$&    & &    &\\
  \hline
  &&&   &     1&&       &\\
  &&&   &     &1&       &\\
  \hline
  &&&   &    & &         1&\\
  &&&   &    & &         &1\\
  \hline
\end{tabular}
$$\\
We will reduce the (1,1) and (2,1) blocks of the
second matrix to the form
\begin{equation}\label{2.2}
\begin{bmatrix}
0&0&0\\0&0&0\\1&0&0
\end{bmatrix}\quad \text{and}\quad
\begin{bmatrix}
0&1&0\\0&0&1
\end{bmatrix},
\end{equation}
preserving the other blocks. Similar to the case
1) of the proof of Theorem \ref{t2.1}, we may add
rows of the (2,1) block to rows of the (1,1)
block.

Suppose $\gamma\ne 0$ (the case $\gamma=0$ and
$\beta\ne 0$ is simpler). Adding the second
column to the first (and making the inverse
transformations with rows to preserve the first
matrix), we make the entry $\beta=0$. Then we
make zero the second column of the (1,1) block.

At last, we interchange the second and third rows
to obtain the blocks \eqref{2.2}. To preserve the
form of the other blocks, we make the same
permutation of columns, then we interchange the
rows and columns within all strips except for the
first horizontal and vertical strips.

Up to simultaneous permutation of rows and
columns, the obtained pair has the form
$\vartriangle^{r+1}\vartriangle^{r+2}$; we have
the bifurcation diagram (ii).\medskip

3) Let ${\cal A}(0,0)=
\vartriangle^r\vartriangle^r \vartriangle^{r+2}$.
Rearranging the rows and columns of ${\cal
A}(\beta,\gamma)$ in the order
$$(\overline{\overline{1}},\overline{\overline{2}}\,|\,
1,\overline{1},\overline{\overline{3}}\,|\,\ldots
\,|\,r, \overline{r},\overline{\overline{r+2}})
\quad\text{and}\quad
(\overline{\overline{1}},\overline{\overline{2}}\,|\,
1,\overline{1},\overline{\overline{3}}\,|\,\ldots
\,|\,r-1,
\overline{r-1},\overline{\overline{r+1}}),$$ we
obtain the following pair of matrices:
$$
\begin{tabular}{|cc|ccc|c|c|}
  \hline
  1&&  &&&    &    \\
  &1&  &&&    &    \\
  \hline
  &&  1&&&   &  \\
  &&  &1&&    &  \\
  &&  &&1&    &  \\
  \hline
  &&  &&&    $\ddots$& \\
  \hline
  &&  &&&    &    $I_3$\\
  \hline
  &&  &&&    &    \\
  \hline
\end{tabular}\qquad
\begin{tabular}{|cc|c|c|c|}
  \hline
  0&0&  & & \\
  1&0&  & &\\
  \hline
  $\beta$&0&  & &\\
  $\gamma$&0&  & &\\
  0&1&  & &\\
  \hline
  &&  $\ddots$&  &\\
  \hline
  &&  & $I_3$ &\\
  \hline
  &&  & &  $I_3$\\
  \hline
\end{tabular}
$$
If $\beta\ne 0$ or $\gamma\ne 0$, then we reduce
the (1,1) and (2,1) blocks of the second matrix
to the form
\begin{equation*}
\begin{bmatrix}
0&0\\0&0
\end{bmatrix}\quad \text{and}\quad
\begin{bmatrix}
0&0\\1&0\\0&1
\end{bmatrix},
\end{equation*}
preserving the other blocks. Up to simultaneous
permutation of rows and columns, the obtained
pair has the form
$\vartriangle^{r}\vartriangle^{r+1}
\vartriangle^{r+1}$; we have the bifurcation
diagram (iii).\medskip

4) Let ${\cal A}(0,0)=
\vartriangle^r\vartriangle^{r+2}
\vartriangle^{r+2}$. We rearrange the rows and
columns of ${\cal A}(\beta,\gamma)$ in the order
$$({\overline{1}},\overline{\overline{1}}\,|\,
\overline{2},\overline{\overline{2}}\,|\,
1,{\overline{3}},\overline{\overline{3}}\,|\, \ldots\,|\,
r,\overline{r+2},\overline{\overline{r+2}})$$ and
$$({\overline{1}},\overline{\overline{1}}\,|\,
\overline{2},\overline{\overline{2}}\,|\,
1,{\overline{3}},\overline{\overline{3}} \,|\,
\ldots\,|\,
r-1,\overline{r+1},\overline{\overline{r+1}})$$
and obtain the following pair of
matrices:\\[-8mm]

$$
\begin{tabular}{|cc|cc|ccc|c|c|}
  \hline
 1&&   &&   &&&   &   \\
 &1&   &&   &&&   &   \\
 \hline
 &&   1&&   &&&   &   \\
 &&   &1&   &&&   &   \\
 \hline
 &&   &&   1&&&   &   \\
 &&   &&   &1&&   &   \\
 &&   &&   &&1&   &   \\
 \hline
 &&   &&   &&&   $\ddots$&   \\
 \hline
 &&   &&   &&&   &   $I_3$\\
 \hline
 &&   &&   &&&   &   \\
 \hline
\end{tabular}\qquad
\begin{tabular}{|cc|cc|c|c|c|}
  \hline
 0&0&   &&   & &   \\
 0&0&   &&   & &   \\
 \hline
 1&&   &&    & &   \\
 &1&   &&    & &   \\
 \hline
 $\beta$&$\gamma$&   &&   & & \\
 &&   1&&   & & \\
 &&   &1&   & & \\
 \hline
 &&   &&   $\ddots$&    &\\
 \hline
 &&   &&   &   $I_3$& \\
 \hline
 &&   &&   &   &   $I_3$\\
 \hline
\end{tabular}
$$
If $\beta\ne 0$ or $\gamma\ne 0$, then we reduce
the (1,1), (2,1), and (3,1) blocks of the second
matrix to the form
\begin{equation*}
\begin{bmatrix}
0&0\\0&0
\end{bmatrix},\quad
\begin{bmatrix}
0&0\\1&0
\end{bmatrix},\quad
\begin{bmatrix}
0&1\\0&0\\ 0&0
\end{bmatrix},
\end{equation*}
preserving the other blocks. We have the
bifurcation diagram (iv).\medskip

5) Let ${\cal A}(0,0)=
\vartriangle^r\vartriangle^{r}\!\!\lambda$.
Suppose $\lambda\ne \infty$ (the case
$\lambda=\infty$ is considered analogously; it
may be also reduced to the considered case by
interchanging the matrices). We rearrange the
rows and columns of ${\cal A}(\beta,\gamma)$ in
the following manner:
$$(\overline{\overline{1}}\,|\, 1,\overline{1}\,|\,
\ldots\,|\, r,\overline{r}) \quad\text{and}\quad
(\overline{\overline{1}}\,|\, 1,\overline{1}\,|\,
\ldots\,|\,
r-1,\overline{r-1}).$$ The obtained pair of
matrices is
$$
\begin{tabular}{|c|cc|c|c|c|}
  \hline
 1&   &&   &   &\\
 \hline
 &   1&&   &   &\\
 &   &1&   &   &\\
 \hline
 &   &&   $I_2$&   &\\
 \hline
 &   &&   &   $\ddots$&\\
 \hline
 &   &&   &   &$I_2$\\
 \hline
 &   &&   &   &\\
 \hline
\end{tabular}\qquad
\begin{tabular}{|c|c|c|c|c|}
  \hline
 $\lambda$&   &   &   &\\
 \hline
 $\beta$&   &   &   &\\
 $\gamma$&   &   &   &\\
 \hline
 &   $I_2$&   &   &\\
 \hline
 &   &    $\ddots$& &\\
 \hline
 &   &    &$I_2$ &\\
 \hline
 &   &   &   &$I_2$\\
 \hline
\end{tabular}\qquad
$$
Let $\beta\ne 0$ or $\gamma\ne 0$. We reduce the
block $[\beta\ \gamma]^T$ to the form $[0\ 1]^T$
by row transformations within the second
horizontal strip (to preserve the form of
matrices, we make the same transformations within
all horizontal strips except for the first strip,
then the inverse transformations within the
vertical strips except for the first strip). At
last, we make the entry $\lambda=0$. We have the
bifurcation diagram (v).\medskip

6) Let ${\cal A}(0,0)=
\vartriangle^r\vartriangle^{r+1}\!\!\lambda$.
Suppose $\lambda\ne \infty$ (the case
$\lambda=\infty$ is considered analogously). We
rearrange the rows and columns of ${\cal
A}(\beta,\gamma)$ in the following manner:
$$(\overline{\overline{1}}\,|\,
\overline{1}\,|\,1,\overline{2}\,|\, \ldots\,|\,
r,\overline{r+1}) \quad\text{and}\quad
(\overline{\overline{1}}\,|\,
\overline{1}\,|\,1,\overline{2}\,|\, \ldots\,|\,
r-1,\overline{r}).$$ The obtained pair of
matrices is
$$
\begin{tabular}{|c|c|cc|c|c|c|}
  \hline
 1&  &  &&  &  &\\
  \hline
 &  1&  &&  &  &\\
  \hline
 &  &  1&&  &  &\\
 &  &  &1&  &  &\\
  \hline
 &  &  &&  $I_2$&  &\\
  \hline
 &  &  &&  &  $\ddots$&\\
  \hline
 &  &  &&  &  &$I_2$\\
  \hline
 &  &  &&  &  &\\
  \hline
\end{tabular}\qquad
\begin{tabular}{|c|c|c|c|c|c|}
  \hline
 $\lambda$&  &  &  &  &\\
  \hline
 $\gamma$&  &  &  &  &\\
  \hline
 $\beta$ &  &  &  &  &\\
 &  1&  &  &  &\\
  \hline
 &  &  $I_2$&  &  &\\
  \hline
 &  &  &  $\ddots$&  &\\
  \hline
 &  &  &  &  $I_2$&\\
  \hline
 &  &  &  &  &$I_2$\\
  \hline
\end{tabular}
$$
If $\beta\ne 0$, then we make the entries
$\beta=1$ and $\gamma=\lambda=0$; the obtained
pair is of type
$\vartriangle^{r+1}\vartriangle^{r+1}$. If
$\beta= 0$ and $\gamma\ne 0$, then we make the
entries $\gamma=1$ and $\lambda=0$; the obtained
pair is of type
$\vartriangle^{r}\vartriangle^{r+2}$. We have the
bifurcation diagram (vi).\medskip

7) Let ${\cal A}(0,0)=
\vartriangle^r\!\!\lambda\mu$. Consider the case
$\lambda\ne \infty$ and $\mu\ne \infty$.
Rearranging the rows and columns of ${\cal
A}(\beta,\gamma)$ in the order
$$(\overline{1},\overline{\overline{1}}\,|\,
1,2,\dots,r)\quad\text{and}\quad
(\overline{1},\overline{\overline{1}}\,|\,
1,2,\dots,r-1),$$ we obtain the following pair of
matrices:
$$
\begin{tabular}{|cc|ccc|}
  \hline
 1&&  &&\\
 &1&  &&\\
  \hline
 &&  1&&\\
 &&  0&$\ddots$&\\
 &&  &$\ddots$&1\\
 &&  &&0\\
  \hline
\end{tabular}\qquad
\begin{tabular}{|cc|ccc|}
  \hline
 $\lambda$&&  &&\\
 &$\mu$&  &&\\
  \hline
 $\beta$&$\gamma$&  0&&\\
 &&  1&$\ddots$&\\
 &&  &$\ddots$&0\\
 &&  &&1\\
  \hline
\end{tabular}
$$

If $\beta\ne 0$ and $\gamma\ne 0$, then we make
the entries $(\beta, \gamma)=(0,1)$ by column
transformations within the first vertical strip
(and by the inverse transformations within the
first horizontal strip to preserve the first
matrix). Since $\lambda\ne\mu$, the (1,1) block
takes the form
$$
\begin{bmatrix}
\lambda&0\\a&\mu
\end{bmatrix}, \quad a\ne 0.
$$
Adding the entry $\gamma$, we make the entry
$\mu=0$; adding $a$, we make $\lambda=0$; then we
make $a=1$. The obtained pair is
$\vartriangle^{r+2}$.

If $\beta= 0$ and $\gamma\ne 0$, then we make
$\gamma=1$ and $\mu=0$ and obtain the pair
$\lambda\! \vartriangle^{r+1}$. If $\beta\ne 0$
and $\gamma= 0$, then the pair is reduced to
$\mu\! \vartriangle^{r+1}$.  We have the
bifurcation diagram (vii).\medskip

8) Let ${\cal A}(0,0)=
\lambda^3_1\lambda_2\ldots\lambda_t$. Consider
the case $\lambda_1\ne\infty$. Then ${\cal
A}(\beta,\gamma)$ has a direct summand $(I_3,A)$,
where
\begin{equation*}\label{2.3}
 A= \begin{bmatrix}
 \lambda_1&1&0\\\beta&\lambda_1&1\\
 \gamma&0&\lambda_1
     \end{bmatrix};
\end{equation*}
we may reduce $A$ by similarity transformations.
The characteristic polynomial of $A$ is $\chi(x)=
x^3-\beta x -\gamma$, the roots of its derivative
$\chi'(x)=3x^2-\beta$ are $\pm\sqrt{\beta/3}$.
The matrix $A$ has multiply eigenvalues if and
only if $\chi(x)$ and $\chi'(x)$ have a common
root; that is,
$$\pm\frac{\beta}{3}
\sqrt{\frac{\beta}{3}}
\mp\beta\sqrt{\frac{\beta}{3}}- \gamma=0,\qquad
\frac{4}{27}\beta^3=\gamma^2.$$ The pair has the
Kronecker type $\mu^2_1\mu_2\ldots\mu_{t+1}$ if
$(4/27)\beta^3=\gamma^2$ and
$\mu_1\mu_2\ldots\mu_{t+2}$ if
$(4/27)\beta^3\ne\gamma^2$. We have the
bifurcation diagram (viii).\medskip

9) Let ${\cal A}(0,0)=
\lambda^2_1\lambda^2_2\lambda_3 \ldots\lambda_t$.
Consider the case $\lambda_1\ne\infty$ and
$\lambda_2\ne\infty$. Then ${\cal
A}(\beta,\gamma)$ has direct summands
$$
\left(\begin{bmatrix}
  1&0\\0&1
\end{bmatrix},
\begin{bmatrix}
  \lambda_1&1\\\beta&\lambda_1
\end{bmatrix}\right)\quad \text{and}\quad
\left(\begin{bmatrix}
  1&0\\0&1
\end{bmatrix},
\begin{bmatrix}
  \lambda_2&1\\\gamma&\lambda_2
\end{bmatrix}\right).
$$
Similar to the case 5) of the proof of Theorem
\ref{t2.1}, the pair has the Kronecker type
$\mu^2_1\mu_2\ldots\mu_{t+1}$ if $\beta=0$ or
$\gamma=0$ and the Kronecker type
$\mu_1\mu_2\ldots\mu_{t+2}$ if $\beta\ne 0$, and
$\gamma\ne 0$. We have the bifurcation diagram
(ix).
\end{proof}

\end{document}